\def\Bbb#1{{\bf #1}}

\def\fnote#1{\footnote}
\def\blacksquare{\hbox{\vrule width 4pt height 4pt depth 0pt}}

\def\cwleftpar#1#2{\leftskip #1 \rightskip #2 plus 1fill}
\def\cwrightpar#1#2{\leftskip #1 plus 1fill \rightskip #2}
\def\cwcenterpar#1#2{\leftskip #1 plus 1fill \rightskip #2 plus 1fill}
\def\cwfullpar#1#2{\leftskip#1\rightskip#2}

\def\cwoutdent#1#2{\llap{\hbox to #1{#2 \hss}}\ignorespaces}
\def\cwparbegin#1#2#3#4#5{
	\ifcase #1 \cwleftpar{#2}{#3}
	\or \cwrightpar{#2}{#3}
	\or \cwcenterpar{#2}{#3}
	\else \cwfullpar{#2}{#3}\fi
	\ifcase #4 \baselineskip = 1.5\baselineskip
	\or \baselineskip = 2\baselineskip
	\or \baselineskip = 3\baselineskip
	\else \baselineskip = 1\baselineskip\fi
	\ifdim #5 > 0in \else \noindent \fi
	\noindent\ignorespaces}
\documentclass{article}
\begin{document}
%------------------------------------------------------------------------
%                               ChiWriter FOOTER
%------------------------------------------------------------------------
\advance \vsize by -1\baselineskip
\def\makefootline{
{\vskip \baselineskip \noindent \folio                                  \par
}}
%------------------------------------------------------------------------
%                               Your Document
%------------------------------------------------------------------------

\vspace*{2ex}

\noindent {\LARGE Flat linear connections in terms
	of\\[1ex] flat linear transports in tensor bundles}

\vspace*{2ex}

\noindent Bozhidar Zakhariev Iliev
\fnote{0}{\noindent $^{\hbox{}}$Permanent address:
Laboratory of Mathematical Modeling in Physics,
Institute for Nuclear Research and \mbox{Nuclear} Energy,
Bulgarian Academy of Sciences,
Boul.\ Tzarigradsko chauss\'ee~72, 1784 Sofia, Bulgaria\\
\indent E-mail address: bozho@inrne.bas.bg\\
\indent URL: http://theo.inrne.bas.bg/$\sim$bozho/}

\vspace*{2ex}

{\bf \noindent Published: Communication JINR, E5-92-544, Dubna, 1992}\\[1ex]
\hphantom{\bf Published: }
http://www.arXiv.org e-Print archive No.~math.DG/0406010\\[2ex]

\noindent
2000 MSC numbers: 53C99, 53B99, 57R35\\
2003 PACS numbers: 02.40.Ma, 02.40.Vh, 04.90.+e\\[2ex]

\noindent
{\small
The \LaTeXe\ source file of this paper was produced by converting a
ChiWriter 3.16 source file into
ChiWriter 4.0 file and then converting the latter file into a
\LaTeX\ 2.09 source file, which was manually edited for correcting numerous
errors and for improving the appearance of the text.  As a result of this
procedure, some errors in the text may exist.
}\\[2ex]

	\begin{abstract}
The parallel linear transports defined by flat linear connection are
axiomatically described. On this basis a number of properties, some of which
are new, of these transports and connections are derived.
	\end{abstract}\vspace{3ex}

 {\bf 1. INTRODUCTION}

\medskip
This work starts investigations devoted to the axiomatic approach to the concept "parallel transport (translation)". In particular, it considers, maybe, the simplest case, namely the one of "flat linear transport over a manifold" in tensor bundles over it which, on the one hand, is sufficiently rich in concrete properties and, on the other hand, admits an "easy" straightforward generalization in different directions.

Section 2 contains the definition of a "flat linear transport" in tensor bundles as a map having the properties described there. This definition is independent of the existence of some additional structures such as metrics or connections. Further, the general form and structure of these transports is derived.

In section 3, it is proved that to any flat linear transport there corresponds a unique flat linear connection defining the parallel transport that coincides with the initial flat linear transport (see propositions $3.1, 3.3$ and 3.6). Moreover, as it is established there, for any flat linear connection there exists a flat linear transport the corresponding to which linear connection coincides with it. Said in other words, this means that the parallel transports generated by flat linear connections are flat linear transports.

Section 4 deals with some local aspects concerning flat linear connections or transports in tensor bundles. Here are derived necessary and sufficient conditions for the existence of local, in some cases holonomic, bases in which the matrix of a given flat linear
transport is constant (and hence unit - see proposition 4.1) or the components of some flat linear connection are zeros. The obtained here results concerning the nonholonomic case are, probably, new as the author failed to find them in the available to him literature.

 In section 5 we present our conclusions.

\medskip
\medskip
 {\bf 2. DEFINITION AND SOME PROPERTIES OF FLAT LINEAR\\
		TRANSPORTS IN TENSOR BUNDLES}

\medskip
Let $M$ be a real smooth, of class $C^{1}$, differentiable manifold [1,2]. By $T^{p,q}_{x}(M)$ we denote the tensor space of type $(p,q)$ over $M$ at $x\in M$; in particular $T^{1,0}_{x}(M)=T_{x}(M)$ and $T^{0,1}_{x}(M)=T^{*}_{x}(M)$ are the tangent and cotangent, respectively, spaces to $M ($see, e.g., [1,2]).

{\bf Definition} ${\bf 2}{\bf .}{\bf 1}{\bf .}$ A flat linear transport (of
tensors) over $M$ is a map $L:(x,y)\mapsto L_{x  \to y}, x,y\in M$, where
$L_{x  \to y}$ is a map from the tensor algebra at $x$ into the tensor
algebra at $y$ having the properties:
\[
L_{x  \to y}(T^{p,q}_{x}(M))\subseteq T^{p,q}_{y}(M),\qquad (2.1)
\]
\[
L_{x  \to y}(\lambda A+\mu A^\prime )=\lambda L_{x  \to y}A+\mu L_{x
\to y}A^\prime ,
\quad  \lambda ,\mu \in {\Bbb R},\ A,A^\prime \in T^{p,q}_{x}(M),
 \qquad (2.2)
\]
\[
L_{x  \to y}(A_{1}\otimes A_{2})
=(L_{x\to y}A_{1})\otimes (L_{x\to y}A_{2}),
\quad  A_{a}\in T^{p_{a}}_{x}(M), a=1,2, \qquad (2.3)
\]
\[
L_{x  \to y}\circ C=C\circ L_{x  \to y},  \qquad (2.4)
\]
\[
L_{y  \to z}\circ L_{x  \to y}=L_{x  \to z},\quad  x,y,z\in M,
 \qquad (2.5)
\]
\[
L_{x  \to x}=id,  \qquad (2.6)
\]
where $C$ is any contraction operator and id means the
identity map (in this case of the tensor algebra at $x)$. The map $L_{x  \to
y}$will be called a flat linear transport from $x$ to y.

{\bf Remark 1.} This definition admits different generalizations to the case of arbitrary fibre bundles but such generalizations will not be needed for the present part of our investigation.

{\bf Remark 2.} As in this work we consider only flat linear transports, we shall call them simply (linear) transports. Here the meaning of the adjective $^{\prime\prime}$flat$^{\prime\prime}$ will be made clear below (see e.g. proposition 3.3).

In other words, we can say that a transport over $M$ is a family of homomorphisms which in fact are isomorphisms (see below), between the tensor algebras at different points of $M$ which preserve the tensor's type, commute with contractions and have the special properties (2.5) and (2.6).

 Putting $z=x$ in (2.5) and taking into account (2.6), we get
 \[
L_{x  \to y  }=L_{y  \to x},\qquad (2.7)
\]
 i.e., the mentioned
homomorphisms have inverse maps which are of the same family, and hence, they
are (linear) isomorphisms.

The following proposition establishes the general functional form of the transports over $M$ as it is specified by $(2.1), (2.2){\bf ,} (2.5)$ and (2.6).

 {\bf Proposition 2.1.} The linear maps $L^{p,q}_{x  \to y}:T^{p,q}_{x}(M)
\to T^{p,q}_{y}(M), x,y\in M$ satisfy (2.5) and (2.6) (with $L^{p,q}_{x  \to
y}$instead of $L_{x  \to y})$ if and only if there exist linear isomorphisms
$L^{p,q}_{x}:T^{p,q}_{x}(M)  \to V, V$ being a vector space, such that
\[
L^{p,q}_{x  \to y}
= \Bigl( L^{p,q}_y \bigr)^{-1} \circ L^{p,q}_{x}.\qquad (2.8)
\]
{\bf
Proof.} Let (2.5) and (2.6) be satisfied by $L^{p,q}_{x  \to y}$. Then, the
substitution of (2.7) into (2.5) gives $L^{p,q}_{x  \to z}=^{  }_{
}L^{p,q}_{z  \to y  }\circ L^{p,q}_{x  \to y}$, for every $x,y,z\in $M.
Therefore, fixing some $x_{0}\in M$, we see that (2.8) is valid for
$V=T^{p,q}_{x_{0}}$and $L^{p,q}_{x}=L^{p,q}_{x  \to x_{0}}$. On the contrary,
if we have the decomposition (2.8), then a straightforward calculation shows
that it converts (2.5) and (2.6) into identities.\blacksquare

{\bf Proposition 2.2.} If the representation (2.8) of $L^{p,q}_{x  \to y}$is
true (see proposition 2.1) and $^\prime V$ is any isomorphic with $V$ vector
space, then
\[
L ^{p,q}_{x  \to }
=
\Bigl( ^\prime L^{p,q}_{y} \Bigr) ^{ -1}\
\circ \Bigl( ^\prime L^{p,q  }_{x}\bigr) ,\qquad (2.9)
\]
 where $^\prime L^{p,q}_{x}:T^{p,q}_{x}(M)  \to ^\prime V$ are isomorphisms, iff  there
exists  an isomorphism $f:V  \to ^\prime V$ such that
\[
^\prime L^{p,q}_{x}=f\circ L^{p,q}_{x}.\qquad (2.10)
\]

 {\bf Proof.} This proposition is almost evident: if (2.10) is true, then
from equation (2.8) it follows (2.9) and vice versa, if (2.8) and (2.9)
hold, then
\(
L^{p,q}_{x  \to y}
=
\Big(  L^{p,q}_y \Big)^{-1} \circ L^{p,q}_{x}
=
\Big( ^\prime L^{p,q}_{y}\Big)^{-1}\circ  ^\prime L^{p,q  }_{x}
\)
and hence
\(
f:=
^\prime L^{p,q  }_{y} \circ \Big( L^{p,q }_{y} \Big)^{-1}
=
^\prime L^{p,q  }_{x}\circ \Big( L^{p,q}_{x}\Big)^{-1}
\)
is the needed isomorphism which does not depend either on $x$ or on
y.\blacksquare

So, if we define $L^{p,q}_{x  \to y}$to be the representation of $L_{x  \to y}$on $T^{p,q}_{x}(M)$, then proposition 2.1 shows that it decomposes according to (2.8) into a composition of two maps depending separately on $y$ and x. The arbitrariness of these last maps is described by proposition 2.2.

Letting $A_{1}=A_{2}=1\in {\Bbb R}$ in (2.3), we find $L_{x  \to y}1=1$
which, by virtue of (2.2), is equivalent to
\[
 L_{x  \to y}\lambda =\lambda , \lambda \in {\Bbb R}.\qquad (2.11)
\]
Let $\{E_{i}(x)\}$ and $\{E^{i}(x)\}$ be dual bases in $T_{x}(M)$ and
$T^{*}_{x}(M)$, respectively, where here and below the Latin indices run from
1 to $n:=\dim(M)$ and the usual summation rule will be assumed. As a
consequence of (2.1) for every $x,y\in M$ there exist uniquely defined
functions $H^{i}_{.j}(y,x)$ and $H^{.j}_{\hbox{i.}}(y,x)$ such that
\[
L_{x  \to y}(E_{j}(x))
= H^{i}_{.j}(y,x)E_{i}(y),
\quad
L_{x \to y}(E^{j}(x)) = H^{.j}_{i.}(y,x)E^{i}(y).
 \qquad (2.12)
\]
If $\delta^{i}_{k}$ are the Kronecker's deltas and $C^{1}_{1}$is the
contraction  operator over the first super- and first subscript, then due to
$(2.4), (2.5)$ and  (2.11),  we  have  $\delta ^{i}_{j}=L_{x  \to y}(\delta
^{i}_{j})=L_{x \to y  }E^{i}(x)(E_{j}(x))=L_{x  \to y
}C^{1}_{1}(E^{i}(x)\otimes E_{j}(x))=C^{1}_{1}\circ L_{x  \to y
}E^{i}(x)\otimes E_{j}(x)=C^{1}_{1}\circ (L_{x  \to
y}E^{i}(x))\otimes (L_{x  \to y}E_{j}(x))=C^{1  }_{1
}(H^{i}_{.k}(y,x)E^{k}(y))\otimes (H^{.l}_{\hbox{j.}}(y,x)E_{l}(y)^{  }_{
}=H^{i}_{.k}(y,x)H^{.k}_{\hbox{j.}}(y,x)$, i.e.
\[
 H^{i}_{.k}(y,x)H^{.k}_{\hbox{j.}}(y,x)=\delta ^{i}_{j}  \qquad (2.13)
\]
 or, using the matrix notation,
\[
H^{i}_{.k}(y,x)  \cdot   H^{.k}_{\hbox{j.}}(y,x)  ={\Bbb I}:=  \delta
^{i}_{j}  ,\qquad (2.13^\prime )
\]
where as a first matrix index is considered
the superscript and as a second one the subscript.

From (2.12) and (2.2) it follows at once that $H^{i}_{.k}(y,x)$ and
$H^{.k}_{\hbox{j.}}(y,x)$ are components of bivectors [6] defined at
$(y,x)\in M  M$, or more precisely, we have
\[
H(y,x):=H^{i}_{.k}(y,x)E_{i}(y)\otimes E^{k}(x)\in T_{y}(M)\otimes
T^{*}_{x}(M),\qquad (2.14a)
\]
\[
H^{-1}(y,x):=H^{.k}_{\hbox{j.}}(y,x)E^{j}(y)\otimes E_{k}(x)\in
T^{*}_{y}(M)\otimes T_{x}(M),\qquad (2.14b)
\]
i.e., $H(y,x)$ is a vector at
$y$ and a covector (1-form) at $x$ and $H^{-1}(y,x)$, its inverse bivector,
is covector at $y$ and vector at x.

The bivectors (2.14) uniquely define the action of $L_{x  \to y}$on any tensor $T\in T^{p,q}_{x}(M)$.

{\bf Proposition 2.3.} If
$T=T^{i_{1}\ldots i_p}_{j_{1}\ldots j_q} E_{i_{1}}(x)\otimes \cdot
\cdot \cdot \otimes E_{i_{p}}(x)\otimes E^{j_{1}}(x)\otimes \cdot \cdot \cdot
\otimes \otimes E^{j_{q}}(x)$, then
\[
L_{x  \to y}(T)
=\Bigl( \prod_{a=1}^{p} H^{k_{a}}_{..i_{a}}(y,x) \Bigr)
 \Bigl(\prod_{b=1}^{q} H^{..j_{b}}_{l_{b}}(y,x) \Bigr)
T^{i_{1}\ldots i_p}_{j_{1}\ldots j_q}
E_{k_{1}}(y)
\]
\[
\otimes \cdot \cdot \cdot \otimes
E_{k_{p}}(y)\otimes E^{l_{1}}(y)\otimes
\cdot \cdot \cdot \otimes E^{l_{q}}(y)
\qquad (2.15)
\]

{\bf Proof.} This result
is a simple corollary from (2.12) and a multiple application of (2.2) and
(2.3).\blacksquare

 If $p=q+1=1$, then from (2.8) and (2.15), we get
\[
L^{1,0}_{x  \to y}(T)=(L^{1,0}_{y})^{-1}\circ
(L^{1,0}_{x})(T)=H^{\hbox{j.}}_{.i}(y,x)T^{i}(x)E_{j}(y).
\]

 Hence, letting $F_{x}:= [(F_{x})^{i}_{.j}] :=[(L^{1,0}_{x})^{i}_{.j}]$ to
be the matrix of the matrix elements of $L^{1,0}_{x}$when some bases
$\{E_{i}(x)\}$ in $T_{x}(M)$ and $\{e_{i}\}$ in $V$ are fixed, i.e.,
$L^{1,0}_{x}(E_{j}(x))=:(L^{1,0}_{x})^{i}_{.j}e_{i}$, and defining $\mathbf{
H}(y,x):=[H^{\hbox{j.}}_{.i}(y,x)]$, we see that
\[
\mathbf{ H}(y,x)=F^{-1}_{y}F_{x},\qquad (2.16)
\]
 where a matrix multiplication is understood.

{\bf Proposition 2.4.} Some map $L_{x  \to y}$of the tensor algebra at $x$ into the tensor algebra at $y$ is a linear transport from $x$ to $y$ if and only if in the corresponding local bases it acts according to (2.15) in which the bivectors (2.14) are inverse to one another, i.e., (2.13) is valid, and (2.16) is true for some nondegenerate matrix $F_{x}$.

{\bf Proof.} If $L_{x  \to y}$is a transport from $x$ to $y$, then, as we already proved, $(2.13)-(2.16)$ are valid, the components of the mentioned bivectors being defined by (2.12), and vice versa, if $(2.13)-(2.16)$ take place, then, as can easily be proved, $(2.1)-(2.6)$ and (2.12) are satisfied for every $F_{x}$, i.e., the so constructed $L_{x  \to y}$is a linear transport from $x$ to y.\blacksquare

 {\bf Proposition 2.5.} Every manifold admits linear transports.

{\bf Proof.} In the proof of proposition 2.3, we saw that to any nondegenerate $n  n$ matrix function $F_{x}$on $M$ and any local basis in its tangent bundle there corresponds, in conformity with (2.15) and $(2.16), a$ linear transport $L_{x  \to y}$from $x$ to $y$ for every $x,y\in $M. So, defining $L:(x,y)\mapsto L_{x  \to y}$, we conclude that $L$ is a linear transport over M.\blacksquare

{\bf Remark.} If $F_{x}$defines some linear transport over $M$, then the
matrix function
\[
^\prime F_{x}=DF_{x}, \det(D)\neq 0,\infty ,\qquad (2.17)
\]
D being a  nondegenerate $n\times n$ constant matrix, defines the same linear
transport, i.$e$ the transport itself defines $F_{x}$up to the constant left
multiplier.  This is a simple corollary from proposition 2.2 (see (2.10)). In
particular, for $^\prime V=V$ the matrix $D$ may be considered as a matrix by
which the basis $\{e_{i}\}$ in $V$ is transformed.

So, as a conclusion of the above discussion, we infer that the definition of a linear transport over $M$ is equivalent to defining in it a pair of inverse to one another bivector fields,  the local representation of which is defined by (2.14) and one of which is given by (2.16).

Below we everywhere assume the manifold $M$ to be endowed with a linear transport L.

\medskip
\medskip
 {\bf 3. THE EQUIVALENCE BETWEEN\\
 FLAT LINEAR  TRANSPORTS IN TENSOR BUNDLES AND\\
 FLAT LINEAR CONNECTIONS}

\medskip
Let us first of all remember some simple facts about linear connections (in tensor bundles) which can be found, e.g., in [1,2].

Let $T^{p,q}(M)$ be the tensor bundle of type $(p,q)$ over M. By Sec$^{k}(T^{p,q}(M))$ and Sec$(T^{p,q}(M))$ we denote, respectively, the set of $C^{k}$and the set of all sections of $T^{p,q}(M)$. Let $T(M)$ be the algebra of tensor fields on M.

From a lot of equivalent definitions of a linear connection on $T(M)$ we choose the following one (see{\bf ,} e.g., [3] or $[2], ch$. III, \S2].

A linear connection on $T(M)$ is a map $\nabla $ such that if $V\in  \in
$Sec$(T^{1,0}(M))$, then $\nabla :V\mapsto \nabla _{V}$where the covariant
derivation (differentiation$) \nabla _{V}$along $V$ has (here by definition)
the properties:

(1) $\nabla _{V}:T(M)  \to T(M)$ is a type preserving derivation, i.e.,
\[
 \nabla _{V}:Sec^{1}(T^{p,q}(M))  \to Sec(T^{p,q}(M)),\qquad (3.1)
\]
\[
  \nabla _{V}\circ C=C\circ \nabla _{V},\qquad (3.2)
\]
\[
 \nabla _{V}(A\otimes B)=(\nabla _{V}A)\otimes B+A\otimes (\nabla
_{V}B),\qquad (3.3)
\]
\[
 \nabla _{V}(A+A^\prime )=\nabla _{V}A+\nabla _{V}A^\prime ,\qquad (3.4)
\]
where $C$ is a contraction operator, $A, B$ and $A^\prime $ are arbitrary
$C^{1}$tensor fields on $M, A$ and $A^\prime $ being of one and the same
type.

(2) If $f:M  \to {\Bbb R}$ is $a C^{1}$function, $V,W\in $Sec$(T^{1,0}(M))$
and A is $a C^{1}$tensor field on $M$, then
\[
 \nabla _{V}f=V(f),\qquad (3.5)
\]
\[
 \nabla _{V+W}=\nabla _{V}+\nabla _{W},\qquad (3.6)
\]
\[
\nabla _{fV}A=f\cdot \nabla _{V}A.  \qquad (3.7)
\]

 If $\{E_{i}\}$ is a field of bases  in a neighborhood of some point of $M$,
then the components (coefficients) $\Gamma ^{j}_{.ki}$of $\nabla $ in it are
defined by
\[
\nabla _{E_{i}}E_{k}=:\Gamma ^{j}_{.ki}E_{j}.\qquad (3.8)
\]
 Every transformation $\{E_{i}\mid _{x}\}  \to \{E_{i^\prime }\mid
_{x}=E^{i}_{i^\prime }(x)E_{i}\mid _{x}\}, x\in M$ leads to the
transformation of $\Gamma ^{i}_{.jk}$into $\Gamma ^{i^\prime }_{..j^\prime
k^\prime }$given by
\[
 \Gamma ^{i^\prime }_{..j^\prime k^\prime }(x)=E^{i^\prime
}_{i}(x)E^{j}_{j^\prime }(x)E^{k}_{k^\prime }(x)\Gamma ^{i}_{.jk}+E^{i^\prime
}_{i}(x)E_{k^\prime }(E^{i}_{j^\prime })\mid _{x},
\qquad (3.9)
\]
 where $[ E^{i^\prime }_{i}(x) ] :=[  E^{i}_{i^\prime }(x)]  ^{-1}$.

Any set of functions $\{\Gamma ^{i}_{.jk}\}$ transforming according to (3.9) defines a unique linear connection whose components in $\{E_{i}\}$ are $\Gamma ^{i}_{.jk}[1,2]$.

 Now we shall turn to the topic of the present section.

Let $L$ be a linear transport over $M, V\in $Sec$(T^{1,0}(M))$ and $S$ be $a C^{1}$tensor field on M.

 We define a map
\[
 \nabla ^{L}:V\mapsto \nabla ^{L}_{V},\qquad (3.10a)
\]
 where $\nabla
^{L}_{V}$maps the $C^{1}$tensor fields on $M$ on the set of tensor fields on
$M$ according to
\[
 (\nabla ^{L}_{V}S)(x)
:=\lim_{\epsilon\to 0 }\Bigl\{ \frac{1}{\varepsilon}
(L_{x_{\epsilon}\to x} S(x_{\epsilon })-S(x)) \Bigr\},\qquad (3.10b)
\]
 where $x\in M$ and in some local
coordinates in a neighborhood of $x$ the coordinates of $x_{\epsilon }$are
$x^{i}_{\epsilon }:=x^{i}+\epsilon V^{i}\mid _{x}$in  which  $\epsilon $
belongs  to  some neighborhood of $0\in {\Bbb R}$ and $V\mid _{x}=V^{i}\mid
_{x}\partial /\partial x^{i}$.

Hereafter, for the existence of the limit in (3.10b) we shall
suppose the transports over $M$ to be smooth, of class $C^{1}$, in a sense that such are the bivectors (2.14) or, equivalently, the matrices $F_{x}, x\in M$ in (2.16).

 From (2.6) and (3.10), we find the simple representation
\[
 (\nabla ^{L}_{V}S)(x)
:= \Bigl[\frac{\partial}{\partial \varepsilon}
  (L_{x_{\epsilon }\to x} S(x_{\epsilon })) \Bigr] \big|_{\epsilon =0}
 \qquad (3.10c)
\]
from where it follows that if $\{x^{i}\}$ are any local coordinates in a
neighborhood of $x$, then the components of $(\nabla ^{L}_{V}S)(x)$ are
\[
 (\nabla ^{L}_{V}S)(x)^{  ...}_{  ...}
= \Bigl[ \frac{\partial}{\partial x^i_\varepsilon}
(L_{x_{\epsilon }\to x}S(x_{\epsilon }))^{...  }_{...  } \Bigr]
\frac{\partial x^i_\varepsilon}{\partial \varepsilon} \Big|_{\epsilon =0}
\]
\[
=V^{i}\mid _{x}\cdot
\Bigl[ \frac{\partial }{\partial y^i}
((L_{y  \to x}S(y))^{...  }_{...}\Bigr]\Big|_{y=x},
 \qquad (3.11)
\]
where by dots we denote the indices corresponding to the type of $S$, and
the components of $L_{y  \to x}S(y)$ are explicitly given by (2.15).

{\bf Proposition 3.1.} The map $\nabla ^{L}$defined by (3.1) is a linear connection on the tensor algebra $T(M)$ which means that $\nabla ^{L}_{V}$satisfies $(3.1)-(3.7)$, i.e., that $\nabla ^{L}_{V}$is a covariant derivation along V.

{\bf Proof.} Eqs. (3.1) and (3.2) are simple corollaries from (3.10b) and (2.1) and (2.4), respectively.

 From (2.2) and (3.10b), we get
\[
\nabla ^{L}_{V}(\lambda A+\mu A^\prime )=\lambda \nabla ^{L}_{V}A+\mu
\nabla ^{L}_{V}A^\prime
 \qquad (3.12)
\]
 which results, for $\lambda =\mu =1$, in (3.4).

The equality (3.3) is a consequence of (2.3) and (3.10b) or also, in a coordinate language, of $(3.11), (2.15)$ and the fact that the components of a tensor product of tensors are the product of the corresponding tensor components.

And at the end, $(3.5)-(3.7)$ follow directly from the local representation (3.11) and, in the case of (3.5), from (2.11).\blacksquare

About the connection $\nabla ^{L}$we shall say that it is generated by
(or associated to) the transport L.

{\bf Proposition 3.2.} If $\{x^{i}\}$ are local coordinates in a
neighborhood of $x\in M$, then in the basis $\{\partial /\partial x^{i}\}$
the components of the linear connection $\nabla ^{L}$ are
\[
 H^{i}_{.jk}(x)
=\frac{\partial H^{i}_{.j}(x),y}{\partial y^k} \Big|_{y=x}
 \qquad (3.13)
\]
 where the matrix $\mathbf{ H}(x,y):=
H^{i}_{.j}(x,y)  $ represents, according to (2.12), the linear transport in
this basis.

{\bf Proof.} If we apply (3.11) to $T\in $Sec$^{1}(T^{1,0}(M))$, then using
(2.15) and (2.16), we get
\[
 [(\nabla ^{L}_{V}T)(x)]^{i}
=V^{k}\mid _{x}  \frac{\partial}{\partial y^k}
[H^{i}_{.j}(x,y)T^{j}(y)]\mid _{y=x}
\]
\[
 =V^{k}\mid _{x} \Bigl[ \frac{\partial}{\partial x^k} T^{i}(x) +
\Bigl(\frac{\partial}{\partial y^k}  (H^{i}_{.j}(x,y) \Bigr)
\Bigr]\Big|_{y=x}T^{j}(x).
\]

So, the comparison of this result with (3.8) shows that in the considered
case the connection's components are exactly (3.13).\blacksquare

{\bf Remark.} If $\{E_{i}\}$ is an arbitrary (local) basis in $T^{1,0}(M)$,
then, as can easily be seen, instead of (3.13), we shall have
\[
H^{i}_{.jk}(x)=[E_{k}\mid _{y}(H^{i}_{.j}(x,y))]\mid _{y=x}.\qquad
(3.13^\prime )
\]
 An important property of a linear connection $\nabla
^{L}$generated by a transport $L$ is that if we put $\mathbf{ H}_{k}(x):=
H^{i}_{.jk}(x)  ^{n}_{i,j=1}$, then due to (2.16) the following local
representation is true
\[
\mathbf{ H}_{k}(x)
= \frac{\partial \mathbf{ H}(x,y)}{\partial y^k}\Big|_{y=x}
=F^{-1}(x)  \frac{\partial F(x)}{\partial y^k}
= - \frac{\partial  \mathbf{ H}(y,x)}{\partial y^k}\Big|_{  y=x},
 \qquad (3.14)
\]

{\bf Proposition 3.3.} The linear connection $\nabla ^{L}$generated by a
linear transport $L$ is flat.

{\bf Proof.} If $\nabla $ is a linear connection, $A,B\in $Sec$(T^{1,0}(M))$
and $[A,B]:=A\circ B-B\circ A$ is the commutator of A and $B$, then the
curvature operator is [2]
\[
R(A,B):=\nabla _{A}\circ \nabla _{B}-\nabla _{B}\circ \nabla _{A}-\nabla
_{[A,B]}.\qquad (3.15)
\]
If $\{E_{i}\}$ is a field of arbitrary local bases,
$[E_{i},E_{j}]=:C^{k}_{ij}E_{k}$and $\Gamma ^{k}_{.ij}$are the components of
$\nabla $ in $\{E_{i}\}$, then the components of the curvature tensor $R$ are
[2,5]
\[
 R^{i}_{\hbox{.jkl}}(x)=-2E_{l}(\Gamma ^{i}_{.jk})|_{x}+\Gamma
^{m}_{.jk}(x)\Gamma ^{i}_{.ml}(x)^{  }_{  [k,l]}-C^{m}_{kl}(x)\Gamma
^{i}_{.jm}(x)  \qquad (3.16)
\]
 where antisymmetrization is performed, e.g.
$(A_{kl})_{[k,l]}:= (A_{kl}-A_{lk})$, over the indices included in
square brackets.

Defining $\Gamma _{k}(x):= [ \Gamma ^{i}_{.jk}(x)]  ^{n}_{i,j=1}$ and
$\mathbf{ R}_{kl}(x):= [ R^{i}_{\hbox{.jkl}}(x)]  ^{n}_{i,j=1}$, where as a first
matrix index is considered the superscript, in any coordinate basis we find
\[
\mathbf{ R}_{kl}(x)
=-2  \Bigl( \frac{\partial\Gamma _{k} }{\partial x^l }
+\Gamma _{l}(x)\Gamma(x)\Bigr)_{ [k,l]}.\qquad (3.16^\prime )
\]
 In particular, for the connection $\nabla ^{L}(3.14)$ is valid the substitution of which
into $(3.16^\prime )$ gives
\[
\mathbf{ R}_{kl}(x)\|_{\nabla =\nabla ^{L}}
=-2 \Bigl[
\frac{\partial }{\partial x^l}
\Bigl( F^{-1}(x)  \frac{\partial F(x) }{\partial x^k} \Bigr)
 -F^{-1}(x) \frac{\partial F(x) }{\partial x^l}
\cdot F^{-1}(x)  \frac{\partial F(x) }{\partial x^k}\Bigr]
\Big|_{[k,l]}\equiv 0,
\]
  where the use of $\partial F^{-1}/\partial x^{k}=-F^{-1}(\partial
F/\partial x^{k})F^{-1}$is made.\blacksquare

{\bf Proposition} ${\bf 3}{\bf .}{\bf 4}{\bf .} A$ linear connection $\nabla $ on $M$ is flat if and only if it is generated by some linear transport $L$, i.e., iff for some transport $L$ we have $\nabla =\nabla ^{L}$.

{\bf Proof.} The sufficiency was already established in proposition 3.3. So,
let's suppose that $\mathbf{ R}_{kl}=0$ for some connection $\nabla $. Then,
there exists a matrix function $F_{x}$such that
\[
 \Gamma _{k}(x)=F^{-1}_{x}\partial F_{x}/\partial x^{k}.\qquad (3.17)
\]
Actually, the integrability conditions for this equation with respect to
$F_{x}$ are
\[
0= \Bigl( \frac{\partial^2 F_x }{\partial x^k\partial x^l}
\Bigr)_{  [k,l]}
=\Bigl[\frac{\partial }{\partial x^k}
( (F_{x}\Gamma _{l}(x))\Bigr]_{ [\hbox{k.}l]}
=-F_{x}\mathbf{ R}_{kl}(x)
\]
 which in the considered here case are satisfied due to $\mathbf{ R}_{kl}=0$.

So, if we define a linear transport $L$ whose matrix in the used basis is $H(y,x):=F^{-1}_{y}F_{x}$, we see that the associated to this transport connection $\nabla ^{L}$has, in accordance with (3.14), components $\mathbf{ H}_{k}(x)=F^{-1}_{x}\partial F_{x}/\partial x^{k}$which by virtue of (3.17) coincide with the ones of $\nabla $, so that $\nabla =\nabla ^{L}.\blacksquare $

In other words, the last proposition states that the definition of a (flat) linear transport in the tensor bundles over $M$ is equivalent to the definition of a flat linear connection in $T(M)$.

 {\bf Proposition 3.5.} If $L$ is a linear transport and
\[
 L_{y}:Sec(T^{p,q}(M))  \to Sec(T^{p,q}(M)),\quad  y\in M\qquad (3.18a)
\]
 is such that for every tensor field A
\[
 (L_{y}A)(x):=L_{y\to x}A(y)), x,y\in M,\qquad (3.18b)
\]
 then
\[
\nabla ^{L}_{V}\circ L_{y}\equiv 0.\qquad (3.19)\]

 {\bf Proof.} This result is a simple corollary from $(3.18b), (2.6)$ and
(3.10).\blacksquare

{\bf Proposition 3.6.} If $L$ is a (flat) linear transport generating the connection $\nabla ^{L}$, then the parallel transport defined by $\nabla ^{L}$coincides with L.

{\bf Proof.} Let $\gamma :J  \to M, J\subset {\Bbb R}$ be $a C^{1}$path and $s,t\in $J.  The  parallel transport for some connection $\nabla $ along $\gamma $ is a map $P^{\gamma }_{\gamma (s),\gamma (t)}:T^{p,q}_{\gamma (s)}  \to    \to T^{p,q}_{\gamma (t)}$such that if $A_{0}\in T^{p,q}_{\gamma (s)}$, then  $P^{\gamma }_{\gamma (s),\gamma (t)}(A_{0})=B_{\gamma (t)}$,  where the tensor field $B$ is defined along $\gamma $ by the initial-value  problem $\nabla _{\cdot }B=0, B_{\gamma (s)}=A_{0}$in which    is the tangent to $\gamma $  vector  field  $(cf. [1-5])$.

As the generated by $L$ connection $\nabla ^{L}$is flat (see proposition
3.3), the defined by it parallel transport does not depend on the path $\gamma $ but only on the points $\gamma (s)$ and $\gamma (t) [4,5]$. This means that the action of this parallel transport is $^{L}P^{\gamma }_{\gamma (s),\gamma (t)}(A_{0})=B_{\gamma (t)}$, where the tensor field $B$ is a solution of $\nabla ^{L}_{V}B=0, B_{\gamma (s)}=A_{0}$for every $V\in $Sec$(T^{1,0}(M))$. By (3.19) and (2.6) this solution is $B=L_{\gamma (s)}A$, where A is any tensor field with the property $A_{\gamma (s)}=A_{0}$. From all this we find $^{L}P^{\gamma }_{\gamma (s),\gamma (t)}(A_{0})=(L_{\gamma (s)}A)(\gamma (t))=L_{\gamma (s)  \to \gamma (t)}A_{0}$and hence, $^{L}P^{\gamma }_{\gamma (s),\gamma (t)}=L_{\gamma (s)  \to \gamma (t)}.\blacksquare $

From propositions 3.4 and 3.6 we infer that any parallel transport defined by a flat linear connection coincides with some flat linear transport and vice versa. This means that the (flat) linear transports in tensor bundles, defined in section 2, realize the axiomatic approach to such parallel transports, i.e., that $(2.1)-(2.6)$, when taken as axioms, define uniquely the set of these parallel transports.

\medskip
\medskip
 {\bf 4. SOME RESULTS CHARACTERIZING THE FLAT CASE}

\medskip
In this section we shall investigate problems concerning the question when in a (local) basis it is possible for the matrix describing in it a (flat) linear transport to be constant or for the components of a (flat) linear connection in it to be zeros.

Elsewhere we shall show that the results presented below are specific of the considered here flat case and that in more general situations they are valid only locally, namely, at a given point or along a given path.

{\bf Proposition 4.1.} For every transport $L$ there exists a field of local bases $\{E_{i^\prime }\}$ in the tangent bundle in which the components of
the corresponding to it bivector $H(x,y)$ are Kronecker's deltas, i.e.,
$H^{i^\prime }_{..j^\prime }(x,y)=\delta ^{i^\prime }_{j^\prime }$. Moreover,
if in some basis the components of $H(x,y)$ are constant (with respect to $x$
and $y)$, then they are Kronecker's deltas and this basis can be obtained from
$\{E_{i^\prime }\}$ through linear transformation with constant coefficients
and on the contrary, in any basis obtained from $\{E_{i^\prime }\}$ by such a
transformation the components of $H(x,y)$ are Kronecker's deltas.

{\bf Proof.} Let $\{E_{i}\}$ be a fixed basis in the tangent to $M$ bundle.
Due to (2.16) there is a matrix function $F_{x}$ such that
\[
H^{i}_{.j}(x,y)
=(F^{-1}_{x})^{i}_{.a}(F_{y})^{a}_{.j}
=\sum_{a=1}^{\dim(M)} (F^{-1}_{x})^{i}_{.a}(F_{y})^{a}_{.j}.
\]

  Let the basis $\{E_{i^\prime }\}$ be defined at any $x\in M$ by
\[
 E_{i^\prime }(x):=\delta ^{a}_{i^\prime
}(F^{-1}_{x})^{i}_{.a}E_{i}(x)\qquad (4.1^\prime )
\]
 Then
\[
 E^{i^\prime }(x)=\delta ^{i^\prime }_{a}(F_{x})^{a}_{.i}E^{i}(x)\qquad
(4.1^{\prime\prime})
\]
 and, because of (2.16), we have
\[
 H^{i^\prime }_{..j^\prime }(x,y)=\delta ^{i^\prime
}_{a}(F_{x})^{a}_{.i}H^{i}_{.j}(x,y)\delta ^{b}_{.j^\prime
}(F^{-1}_{y})^{j}_{.b}=\delta ^{i^\prime }_{a}\delta ^{a}_{j^\prime }=\delta
^{i^\prime }_{j^\prime },
\]
i.e., $\{E_{i^\prime }\}$ is a basis with the needed properties.

Let $\{E_{i}\}$ be a fixed basis in which the transport is described by the
matrix $\mathbf{ H}(x,y)$. Then, there exists a nondegenerate matrix $A(x)=
A^{i^\prime }_{i}(x)  :=  A^{j}_{j^\prime }(x)  ^{-1}$such that
$E_{i}(x)=A^{i^\prime }_{i}(x)E_{i^\prime }(x)$ and $E^{i^\prime }(x)=
=A^{i^\prime }_{i}(x)E^{i}(x)$. As a consequence of this, (2.12) and the
above definition of $E_{i^\prime }(x)$, we have $\mathbf{ H}(x,y)=A(x)
H^{i^\prime }_{..j^\prime }(x,y)  (A(y))^{-1}=A(x){\Bbb
I}(A(y))^{-1}=A(x)(A(y))^{-1}$. From here it immediately follows that $\mathbf{
H}(x,y)=$const if and only if $A(x)=$const and if this is the case, then,
evidently, we have $\mathbf{ H}(x,y)={\Bbb I}.\blacksquare $

{\bf Remark.} A basis $\{\tilde{E}_{i^\prime }\}$ with the property
described in proposition 4.1 can be constructed also in the following way.
Take any fixed field $\{E_{i}\}$ of local bases, fix a point $z\in M$ and
define (see (3.18))
\[
\tilde{E}_{i^\prime }
:=\delta ^{i}_{i^\prime }L_{z}(E_{i}).
\qquad (4.1^{\prime\prime\prime})
\]
 Then, due to (2.5) and (3.18a), we have
$\tilde{E}_{i^\prime }(y)=L_{x  \to y}\tilde{E}_{i^\prime }(x)$, hence $\mathbf{
H}^\prime (x,y)={\Bbb I}$. The second part of the proposition can also be
easily derived from $(4.1^{\prime\prime\prime})$.

In the general case the bases in which the components of $H(x,y)$ are
constant are nonholonomic, i.e., they are not generated by some local
coordinates [6]. In a formal language this is expressed by

{\bf Proposition  4.2.} In the tangent bundle $T(M)$ there exists a field of
local holonomic bases, i.e., bases generated by some local coordinates in
$M$, in which the components of the bivector $H(x,y)$ are constant if and
only if in $M$ there exist local coordinates $\{x^{i}\}$ such that in the
associated to them basis $\{\partial /\partial x^{i}\}$ the components of
one, and hence of all, matrix $F_{x}=  (F_{x})^{i}_{.j}  $, defining through
(2.16) the transport $L$ in it, satisfy the equations
\[
\bigl(
\partial (F_{x})^{i}_{.j}/\partial x^{k  }
\big)_{[j,k]}=0.\qquad (4.2)
\]
 Moreover, if the described above coordinates exist, then any local basis,
in which the components of $H(x,y)$ are constant, is holonomic.

 {\bf Remark.} If we define the 1-forms
\[
 \tilde{F}^{i}_{x}:=(F_{x})^{i}_{.j}E^{j}(x),
\tilde{H}^{i}(y,x):=H^{i}_{.k}(y,x)E^{k}(x)
=(F^{-1}_{y})^{i}_{.j}(\tilde{F}^{j}_{x})
\]
then, as can easily be proved, (4.2) is equivalent to the statement that
any one of these forms is closed, i.e., to either of
\[
d\tilde{F}^{i}_{x}=0,\qquad (4.3)
\]
\[
d_{x}(\tilde{H}^{i}(y,x))=0,\qquad (4.3^\prime )
\]
where $d_{x}$means exterior derivation with respect to x.

{\bf Proof.} Let us take the basis $\{\partial /\partial x^{i}\}$ associated
with some fixed local coordinates $\{x^{i}\}$. If in $(4.1^\prime )$ and
$(4.1^{\prime\prime})$ we substitute $E_{i}(x)=\partial /\partial x^{i}$,
then, by proposition 4.1 and its proof, any basis in which the components of
$H(x,y)$ are constant is of the form
\[
 \tilde{E}_{j}(x)
=A^{i^\prime }_{j}\delta ^{a}_{i^\prime }
(F^{-1}_{x})^{i}_{.a}(\partial /\partial x^{i}),\qquad (4.4^\prime )
\]
\[
\tilde{E}^{j}(x)
=A^{j}_{i^\prime }\delta ^{i^\prime }_{a}(F_{x})^{a}_{.i}dx^{i},
\qquad (4.4^{\prime\prime})
\]
where $[ A^{j}_{i^\prime }] =[A^{i^\prime }_{j}]^{-1}$ is a constant matrix
and the components of $F_{x}$and $F^{-1}_{x}$ are referred to
$\{\partial /\partial x^{i}\}$ and $\{dx^{i}\}$.

By definition the bases (4.4) are holonomic if the 1-forms
$(4.4^{\prime\prime})$ are exact [6], i.e., if there exist
$\tilde{x}^{j}=\tilde{x}^{j}(x)$, such that
$\tilde{E}^{j}(x)=d\tilde{x}^{j}(x)$, or equivalently $\tilde{E}_{j}=\partial
/\partial \tilde{x}^{j}$, i.e., that $\{\tilde{x}^{j}\}$ may be taken as
local coordinates. Locally, a necessary and sufficient condition for that is
d\~{E}$_{j}(x)=0 ($see the converse of Poincar\`e's lemma in [3], p. $145;
cf. [1,6])$ which, as can easily be seen by means of $(4.4^{\prime\prime})$,
is equivalent to (4.3).

This proves the first part of the proposition. Its second part is a trivial corollary from the second part of proposition 4.1 and the evident fact that a linear combination with constant coefficients of exact 1-forms is an exact 1-form.\blacksquare

Before the formulation of the next proposition, which shows the meaning of
proposition 4.1 in terms of connections generated by a linear transports, we
would like to remind $(cf. [1])$ that the torsion tensor $T$ of a linear
connection $\nabla $ is defined by
\[
 T(A,B):=\nabla _{A}B-\nabla _{B}A-[A,B],\qquad (4.5)
\]
$A$ and $B$ being vector fields, and in a local basis $\{E_{i}\}$ its
components are
\[
 T^{i}_{.jk}=-2(\Gamma ^{i}_{.jk})_{[j,k]}-C^{i}_{.jk}.\qquad (4.6)
\]

 {\bf Proposition 4.3.} The torsion $^{L}T$ of the connection $\nabla
^{L}$associated with some linear transport $L$ vanishes if and only if the
conditions (4.2) are fulfilled.

{\bf Proof.} As in a local coordinate basis the connection coefficients of
$\nabla ^{L}$are (see (3.14))
\[
  ^{L}\Gamma ^{i}_{.jk}(x)=H^{i}_{.jk}(x)=(F^{-1}_{x})^{i}_{.a}(\partial
(F_{x})^{a}_{.j}/\partial x^{k}),\qquad (4.7)
\]
 the local components of $^{L}T$ in the same basis are
\[
^{L}T^{i}_{.jk}(x)=-(^{L}\Gamma ^{i}_{.jk}(x))_{[j,k]}=(F^{-1}_{x})^{i}_{.a
}\partial (F_{x})^{a}_{.j}/\partial x^{k  }.\qquad (4.8)
\]
 The comparison of
this result with (4.2) shows the equivalence of (4.2) and the equality
$^{L}T=0.\blacksquare $

{\bf Corollary 4.1.} The torsion $^{L}T$ of the associated with a linear transport $L$ connection $\nabla ^{L}$to be zero is a necessary and sufficient condition for the existence of a field of local holonomic bases in which the matrix (2.16), describing that transport $L$, is constant.

 {\bf Proof.} This result follows from propositions 4.2 and 4.3.\blacksquare

{\bf Corollary 4.2.} If the curvature of a linear connection is zero, then a necessary and sufficient condition for the existence of a local holonomic basis in which the components of the connection are zeros is its torsion to vanish.

{\bf Remark 1.} If the connection is not curvature free, then due to $(3.16) a$ basis with the described property does not exist (see also below corollary 4.3).

{\bf Remark 2.} This is an old classical result which in a somewhat different formulation can be found, for instance, in [6], p. 142 or
in $[5], \S106$, p. 519.

{\bf Proof.} As for a holonomic basis $C^{i}_{.jk}=0$, the necessity
directly follows from (4.6).

On the contrary, let $T^{i}_{.jk}=0$. As $R=0$, then by proposition 3.4
there exists a transport $L$ such that $\nabla =\nabla ^{L}$. But then
$^{L}T^{i}_{.jk}=T^{i}_{.jk}=0$ and due to proposition 4.3 there is a local
holonomic basis in which the components $H^{i}_{.j}(y,x)$ of the matrix
representing the transport in it are constant. So, due to (3.13), in this
local holonomic basis $\Gamma ^{i}_{.jk}(x)=^{L}\Gamma
^{i}_{.jk}(x)=H^{i}_{.jk}(x)=\partial H^{i}_{.j}(x)/\partial
x^{k}=0.\blacksquare $

{\bf Corollary} ${\bf 4}{\bf .}{\bf 3}{\bf .} A$ linear connection is curvature free if and only if there exists a basis in which its components are zeros.

{\bf Remark 1.} Corollary 4.2 tells us when the mentioned basis is holonomic.

{\bf Remark 2.} If one considers only holonomic coordinates [6], then this is a known result; $cf. [6]$, p.142 or $[5], \S106$, p.519.

{\bf Proof.} If the connection $\nabla $ is curvature free, then by
proposition 3.4 there is a transport $L$ such that $\nabla =\nabla ^{L}$. For
$L$, by proposition 4.1, there exists a basis $\{E_{i}\}$ in which the
defined through it matrix (2.16) is constant. In this basis, by $(3.13^\prime
)$ the components of $\nabla $ are $\Gamma
^{i}_{.jk}(x)=H^{i}_{.jk}(x)=E_{k}(y)(H^{i}_{.j}(x,y))=0$, as in it
$H^{i}_{.j}(x,y)=$const.

On the contrary, if for $\nabla $ there exists a basis in which its
components are zeros, then from (3.6) it follow that $\nabla $ is curvature
free.\blacksquare

\medskip
\medskip
 {\bf 5. COMMENTS.}

\medskip
In this work, we have axiomatically defined $^{\prime\prime}$flat linear transports$^{\prime\prime}$ in tensor bundles the class of which, as was proved, coincides with the one of parallel transports generated by flat linear connections. A feature of our approach is that we have fixed for the mentioned definition only those properties of the latter transports which describe them completely. This consideration of parallel transports generated by flat linear connections turns out to be rather fruitful because it is independent of the standard connection theory and it gives possibilities for different generalizations which will be a subject for forthcoming papers.

On the basis of the developed formalism we have expressed a number of properties of flat linear connections in terms of flat linear transports. As the latter are global (integral) objects, the proofs of these properties are considerably simplified with respect to the ones made by means of connections.

Possibly, part of the mentioned properties are new at least in their formulation, but the proofs of all of them are new in spite of that some of them are close to the ones in the references.

\medskip
\medskip
 {\bf ACKNOWLEDGEMENT}

\medskip
This research is partially supported by the Foundation for Scientific Research of Bulgaria under contract $F 103$.

\medskip
\medskip
 {\bf REFERENCES}

\medskip
1. Dubrovin B., S. P. Novikov, A. Fomenko, Modern Geometry, I. Methods and Applications, Springer Verlag.\par
2. Kobayashi S., K. Nomizu, Foundations of Differential Geometry, Vol. 1, Interscience Publishers, New York-London, 1963.\par
3. Lovelock D., H. Rund, Tensors, Differential Forms, and  Variational Principals, Wiley-Interscience Publication, John Wiley \& Sons, New York-London-Sydney-Toronto, 1975.\par
4. Norden A. P., Spaces with Affine Connection, Nauka, Moscow, 1976 (In Russian).\par
5. Rashevskii P. K., Riemannian Geometry and Tensor Analysis, Nauka, Moscow, 1967 (In Russian).\par
6. Schouten J. A., Ricci-Calculus: An Introduction to Tensor Analysis and its Geometrical Applications, $2-nd ed$., Springer Verlag, Berlin-G\"ottingen-Heidelberg, 1954.
\par

\newpage
\vspace{8ex}

\noindent
 Iliev B. Z.\\[5ex]

\noindent
 Flat linear connections in terms of\\
 flat linear transports  in tensor bundles \\[5ex]

\medskip
The parallel linear transports defined by flat linear connection are
axiomatically described. On this basis a number of properties, some of which
are new, of these transports and connections are derived.\\[5ex]

The investigation has been performed at the Laboratory for Theoretical
Physics, JINR.

\end{document}